\documentclass[12pt]{amsart}

\newsymbol\pp 1275
\newsymbol\twoheadrightarrow 1310
\newcommand{\Tor}{\operatorname{Tor}}

\newtheorem{theorem}{Theorem}
\newtheorem{proposition}[theorem]{Proposition}

\newtheorem{conjecture}[theorem]{Conjecture}

\theoremstyle{definition}

\newtheorem{example}[theorem]{Example}

\newcommand{\reg}{\operatorname{reg}}

\title{Degree Bounds for Syzygies of Invariants}
\author{Harm Derksen}
\thanks{The author is supported by NSF grant DMS 0102193.}
\begin{document}
\maketitle
\begin{abstract}
Suppose that $G$ is a linearly reductive group.
Good degree bounds for generators of invariant
rings  were given in \cite{Derksen2}.
Here we study  minimal free resolutions of invariant rings.
For finite linearly reductive groups $G$
it was recently shown in~\cite{Fleischmann,Fogarty,DerksenSidman}
that  rings of invariants
are generated in degree at most the group order $|G|$. 
In characteristic 0 this degree bound is a classical result by Emmy Noether
 (see~\cite{Noether}).
 Given an  invariant ring
 of a finite linearly reductive group $G$,
 we prove that
 the ideal of relations of a minimal set of
generators is generated in degree at most $\leq 2|G|$
\end{abstract}

\section{Introduction}
Let us fix a linearly
reductive group $G$ over a field $K$. If $V$ is a representation
of $G$, then $G$ acts on the coordinate ring $K[V]$. The ring of
invariant functions is denoted by $K[V]^G$.
Let $\beta_G(V)$ be the smallest positive integer $d$
such that all invariants of degree at most $d$ generate the invariant ring.

Noether proved that $\beta_G(V)\leq |G|$
for a finite groups $G$ in characteristic $0$.
A group $G$ over a field $K$ is linearly reductive if and only
if the characteristic of the base field does not divide $|G|$.
This situation is often referred to as the {\it non-modular case}.
Noether's bound  was recently extended
by Fleischmann (see~\cite{Fleischmann}) to the general non-modular case.
Fogarty found another proof of this independently in \cite{Fogarty}. A third
proof follows from the subspaces conjecture in~\cite{Derksen}, which
was solved by  Sidman and the author in~\cite{DerksenSidman}.
For connected linearly reductive groups, 
upper bounds for $\beta_G(V)$ which depend polynomially
on the weights appearing in the representation were given in~\cite{Derksen2}.

In this paper, we discuss good degree bounds for  syzygies of
 invariant rings. Suppose that $G$ is a linearly reductive group
 and that $V$ is a $n$-dimensional representation.
Let $f_1,f_2,\dots,f_r$ be minimal generators of the
invariant ring $R=K[V]^G$ whose degrees
are $d_1\geq d_2\geq \cdots \geq d_r$. Define the graded polynomial
ring $S=K[x_1,x_2,\dots,x_r]$ by $\deg(x_i)=d_i$ for all $i$.
Now $R$ is an $S$-module via the surjective ring homomorphism 
$\varphi:S\twoheadrightarrow R$ defined by $\varphi(x_i)=f_i$ for all $i$.
Let 
$$
0\to F_k\to F_{k-1}\to \cdots \to F_1\to F_0\to R\to 0
$$
be the minimal free graded resolution of $R$ as an $S$-module. 
We define $\beta_G^i(V)$ as the smallest integer such that
$F_i$ is generated as an $S$-module in degree at most $d$.
Note that $F_0=S$ and the image of $F_1$ in $F_0$ is the syzygy ideal
$J$ of $S$ defined by
$$
J=\{h\in k[x_1,\dots,x_r]\mid h(f_1,f_2,\dots,f_r)=0\}.
$$
The invariant ring $K[V]^G$ is Cohen-Macaulay (see \cite{HochsterRoberts}).
By standard methods, one can estimate the degrees in the minimal
resolution of $K[V]^G$ by cutting the invariant ring down
by hypersurfaces. Using an estimate of Knop for
the $a$-invariant (see~\cite[Definition 4.3.6, p. 169]{BH} for
a discussion of the $a$-invariant)
 of the invariant ring, we will prove the following bound:
\begin{theorem}\label{theo1}
We have
$$
\beta_G^i(V)\leq d_1+d_2+\cdots+d_{s+i}-s\leq (s+i)\beta_G(V)-s
$$
where $s$ is the Krull dimension of $K[V]^G$.
\end{theorem}
If $G$ is connected and linearly reductive,
then Theorem~\ref{theo1} and
the polynomial bounds of $\beta_G(V)$
in \cite{Derksen2} give polynomial bounds for $\beta_G^i(V)$.
From the case $i=1$ follows that the syzygy ideal
 $J$ is generated in degree at most $\beta_G^1(V)\leq (s+1)\beta_G(V)-s$.
 
 Let $G$ be a finite group in the nonmodular case.
 Theorem~\ref{theo1} and the inequality $\beta_G(V)\leq |G|$ imply that
 $$\beta_G^i(V)\leq (n+i)|G|-n$$
 and that $J$ is generated in degree
 at most $\beta_G^1(V)\leq (n+1)|G|-n$. 
 This last inequality will be improved in Theorem~\ref{theo2}.
 
 Suppose that $M$ is a graded module over the
 ring $T=K[V]$ with minimal free resolution
 $$
 0\to H_l\to H_{l-1}\to \cdots \to H_0\to M\to 0.
 $$
 Recall that the Castelnuovo-Mumford regularity $\reg(M)$ is the smallest
 integer $d$ such that  $H_i$ is generated in degree at most $d+i$
 for all $i$.
 Let $I$ be the ideal generated by all homogeneous invariants
 of positive degree. Define
 $\tau_G(V)$ as the smallest  integer $d$ such that
 every homogeneous polynomial of degree $d$ lies in  $I$.
 It is well-known that the Castelnuovo-Mumford regularity
 of a finite length graded $T$-module $M$ is exactly the maximum
 degree appearing in $M$ (see~\cite[Exercise 20.15]{Eisenbud}
 or Theorem~\ref{theo6}). 
Application to the module $T/I$ gives us $\reg(T/I)=\tau_G(V)-1$. 
 From \cite[Corollary 20.19]{Eisenbud}
 and the exact sequence $0\to I\to T\to T/I\to 0$
 follows that
 $$
 \reg(I)=\reg(T/I)+1=\tau_G(V).
 $$
  Fogarty's proof of the Noether bound 
 in the nonmodular case (see~\cite{Fogarty}),
  shows that
 $\tau_G(V)\leq |G|$. The following theorem gives a similar
 bound for the degrees of the generators of the syzygy ideal.
 \begin{theorem}\label{theo2}
Suppose that $G$ is a finite group in the nonmodular case. 
Suppose that $\{f_1,f_2,\dots,f_r\}$ is a  minimal set of
homogeneous
generators of the invariant ring $K[V]^G$ and let
$J\subseteq K[x_1,x_2,\dots,x_r]$ be the syzygy ideal.
Then $J$ is generated in degree at most
$$
2\tau_G(V)\leq 2|G|.
$$
\end{theorem}
The proof of this theorem does not seem to extend to higher syzygies.
We conjecture though that similar bounds will hold for higher syzygies:
\begin{conjecture}
If $G$ is a finite group in the nonmodular case, then
$$
\beta^i_G(V)\leq (i+1)\tau_G(V)\leq (i+1)|G|.
$$
\end{conjecture}
\section{Examples}
\begin{example}
Let $S_n$ be the symmetric group acting on $V=K^n$ by
permutation of the coordinates
and let $A_n\subset S_n$ be the alternating group.
The coordinate ring $K[V]$ can be identified with
the polynomial ring $K[y_1,\dots,y_n]$.
The invariant ring of $S_n$ is 
$$K[y_1,y_2,\dots,y_n]^{S_n}=K[e_1,e_2,\dots,e_n]
$$
where $e_j$ is the elementary symmetric polynomial
of degree $j$ defined by
$$
e_j=\sum_{1\leq i_1<i_2<\cdots<i_j\leq n}
y_{i_1}y_{i_2}\cdots y_{i_j}
$$
for all $j$.

It is also well known that the ring of $A_n$-invariants
is
$$
K[y_1,\dots,y_n]^{A_n}=K[e_1,e_2,\dots,e_n,\Delta]
$$
where $\Delta$ is the $A_n$-invariant of degree $n(n-1)/2$ defined by
$$
\Delta=\prod_{1\leq i<j\leq n}(y_i-y_j).
$$
Let us again define a surjective ring homomorphism
$$
K[x_1,x_2,\dots,x_{n+1}]\to K[y_1,\dots,y_n]^{A_n}
$$
where $x_i$ maps to $e_i$ for $i\leq n$ and
$x_{n+1}$ maps to $\Delta$.
The kernel of the homomorphism is the syzygy ideal $J$.

Theorem~\ref{theo1} says that $J$ is generated in degree at most
$$
n(n-1)/2+n+(n-1)+\cdots+1-n=n(n-1).
$$

The polynomials $e_1,e_2,\dots,e_n$ form a regular sequence
in $K[y_1,\dots,y_n]$ and the Hilbert series of
$K[y_1,\dots,y_n]/(e_1,e_2,\dots,e_n)$ is
$$
(1+t)(1+t+t^2)\cdots (1+t+\cdots+t^{n-1})=1+(n-1)t+\cdots+t^{n(n-1)/2}.
$$
From this we see that the highest
degree for which  $K[y_1,\dots,y_n]/(e_1,e_2,\dots,e_n)$
is nonzero is $n(n-1)/2$. The degree $n(n-1)/2$ part
of $K[y_1,\dots,y_n]/(e_1,\dots,e_n)$ is one-dimensional and
 and it is spanned by the invariant $\Delta$.
This shows that every polynomial of degree at least $n(n-1)/2$
lies in the ideal $(e_1,e_2,\dots,e_n,\Delta)$, so
$\tau_{A_n}(V)=n(n-1)/2$.
By Theorem~\ref{theo2} we get that $J$ is generated in degree at most
$$
 2\cdot n(n-1)/2=n(n-1).
$$
Both theorems are sharp in this example. We have that $\Delta^2$
is $S_n$-invariant and therefore $\Delta^2$ is a polynomial
in $e_1,e_2,\dots,e_n$. This gives a relation of degree $n(n-1)$
and it is known that this relation generates the ideal $J$.
\end{example}
\begin{example}
Let $G$ be the cyclic group of order $m$, generated by $\sigma$.
Let $\sigma$ act on $V=K^n$ by scalar multiplication
with a primitive $m$-th root of unity $\zeta$. This defines
a group action of $G$ on $V$. We identify again $K[V]=K[y_1,\dots,y_n]$
where $y_i$ is the $i$-th coordinate function.
The invariant ring
$K[V]^G$ is generated by the set ${\mathcal M}$
of all monomials in $y_1,y_2,\dots,y_n$
of degree $m$. To every monomial $M\in {\mathcal M}$  we
attach a formal variable $x_M$. We consider the surjective
ring homomorphism
$$
K[\{x_M\}_{M\in {\mathcal M}}]\twoheadrightarrow K[V]^G
$$
which maps $x_M$ to $M$ for every monomial $M\in {\mathcal M}$.
The kernel of this homomorphism is again the syzygy ideal $J$.

By Theorem~\ref{theo1}, $J$ is generated in degree at most $(n+1)m-n$.
Since every monomial of degree $m$ lies the ideal $J$ generated
by all homogeneous invariants of positive degree,
we have $\tau_G(V)=m$.
By Theorem~\ref{theo2}, $J$ is generated in degree at most $ 2m$
(which means that $J$ is generated by polynomials which
are quadratic in the variables $\{x_M\}_{M\in {\mathcal M}}$).
Indeed, $J$ is generated by relations of the form
$$
x_{y_iM} x_{y_jN}-x_{y_jM}x_{y_iN}
$$
where $M$ and $N$ are monomials of degree $m-1$. Now Theorem~\ref{theo2}
is sharp, but Theorem~\ref{theo1} is not.
\end{example}

\section{A general degree bound for Syzygies}\label{sec3}
Suppose that $S=K[x_1,x_2,\dots,x_r]$
is the graded polynomial ring where $\deg(x_i)=d_i$
is a positive integer
for all $i$. We will assume that $d_1\geq d_2\geq \cdots \geq d_r$.
Let $M$ be a finitely generated graded Cohen-Macaulay $S$-module.
The minimal resolution of $M$ is 
$$
0\to F_k\to F_{k-1}\to \cdots \to F_1\to F_0\to M\to 0
$$
where
$$
F_i\cong\Tor^S_i(M,K)\otimes_K S.
$$
Here $\Tor^S_i(M,K)$ is a finite dimensional graded vector space,
and this makes \break $\Tor^S_i(M,K)\otimes_K S$ into a graded module.

 If $M$ is a finite dimensional
graded vector space, then $\deg(M)$ is the maximal degree appearing
in $M$ if $M$ is nonzero, and $\deg(M)=-\infty$ if
$M$ is zero. For a finitely generated module $M$, 
$a(M)$ is the degree of the Hilbert
series $H(M,t)$, seen as a rational function (the so-called $a$-invariant of
$M$).
\begin{theorem}\label{theo6}
We have the inequality
$$
\deg(\Tor^S_i(M,K))\leq d_1+d_2+\cdots+d_{s+i}+a(M)
$$
where $s$ is the dimension of $M$.
\end{theorem}
\begin{proof}
We prove the theorem by induction on $s=\dim M$.
Suppose that $M$ has dimension 0. In this case we prove the
inequality by induction of the length $\dim_K M$ of $M$.
If $M$ has length 0, then $M$ is the trivial module and the
inequality is obvious. Suppose that $M$ is nonzero.
 Note that $a:=a(M)$ is the
maximum degree appearing in $M$. 
Let $M_a$ be the part of $M$ of degree $a$.
Then $M_a$ is a submodule of $M$.
We have an exact sequence of $S$-modules
$$
0\to M_a\to M\to M/M_a\to 0
$$
Since $\dim_K M/M_a<\dim_K M$ and $a(M/M_a)<a$ we get by induction that
$$
\deg(\Tor^S_i(M/M_a,K))\leq d_1+\cdots+d_i+a-1.
$$
The submodule $M_a$ is isomorphic to the module $K^m[-a]$
which is the module $K^m$ whose degree is shifted by $a$.
Since
$$
\deg(\Tor^S_i(K,K))\leq d_1+d_2+\cdots+d_i
$$
by the Koszul resolution, we have that
$$
\deg(\Tor^S_i(M_a,K))\leq d_1+d_2+\cdots+d_i+a.
$$
From the long exact sequence
$$
\cdots \to \Tor^S_i(M_a,K)\to \Tor^S_i(M,K)\to \Tor^S_i(M/M_a,K)\to \cdots
$$
follows that
$$
\deg(\Tor^S_i(M,K))\leq d_1+d_2+\cdots+d_i+a.
$$

Now suppose that $s>0$. Since $M$ is Cohen-Macaulay
we can find a homogeneous nonzero divisor $p$ of
degree $e>0$
and $M/pM$ is again Cohen Macaulay.
First, note that $H(M/pM,t)=(1-t^e)H(M,t)$,
so $a(M/pM)=a(M)+e$.
From the short exact sequence
$$
0\to M[-e]\to M\to M/pM\to 0
$$
we obtain a long exact sequence
$$
\cdots\to\Tor^S_{i+1}(M/pM,K))\to \Tor^S_i(M,K)[-e]\to \Tor^S_i(M,K)\to\cdots
$$
Any element of $\Tor^S_i(M,K)[-e]$ of maximal degree must map to 0 in
$\Tor^S_i(M,K)$, and therefore it must come
from $\Tor^S_{i+1}(M/pM,K)$.
This shows that
$$
e+\deg(\Tor^S_i(M,K))=\deg(\Tor^S_i(M,K)[-e])\leq
\deg(\Tor^S_{i+1}(M/pM,K))\leq
$$
$$
\leq d_1+d_2+\cdots+d_{(s-1)+(i+1)}+a(M/pM)=
d_1+d_2+\cdots+d_{s+i}+a(M)+e,
$$
so finally
$$
\deg(\Tor^S_i(M,K))\leq d_1+d_2+\cdots+d_{s+i}+a(M).
$$
\end{proof}
\begin{proof}[Proof of Theorem~\ref{theo1}]
Let us choose $M=R$ in the previous theorem.
Then 
$$
\beta^i_G(V)=\deg(\Tor^S_i(M,K))\leq d_1+d_2+\cdots+d_{s+i}+a(R)
$$
Knop proved that $a(R)\leq -s$ (see~\cite[Satz 4]{Knop}) 
and Theorem~\ref{theo1} follows.
\end{proof}
\section{bounds for the syzygy ideal for finite groups}
\begin{proposition}
Suppose that  $R=\bigoplus_{d\geq 0} R_d$ is a graded ring
with $R_0=K$
and that $\{f_1,f_2,\dots,f_r\}$ is a minimal
set of homogeneous generators of $R$.
Let $S=K[x_1,\dots,x_r]$ be the
graded polynomial ring and let $\varphi:S\twoheadrightarrow R$
be the surjective ring homomorphism defined by $x_i\mapsto f_i$
for all $i$.
We have an exact sequence of graded vector spaces
$$
\Tor^S_2(K,K)\to \Tor_2^R(K,K)\to \Tor_1^S(R,K)\to 0
$$
\end{proposition}

\begin{proof}
From Exercise A3.47 (with the role of $R$ and $S$ interchanged) 
in \cite{Eisenbud}, we get a five-term exact sequence
$$
\Tor^S_2(K,K)\to \Tor_2^R(K,K)\to \Tor_1^S(R,K)\to \Tor^S_1(K,K)\to
\Tor^R_1(K,K)\to 0.
$$
Let ${\mathfrak n}=(x_1,\dots,x_r)$ be the maximal homogeneous
ideal of $S$ and let ${\mathfrak m}=(f_1,\dots,f_r)$
be the maximal homogeneous ideal of $R$.
Now $\Tor^S_1(K,K)$ and $\Tor^R_1(K,K)$ can be identified
with ${\mathfrak n}/{\mathfrak n}^2$ and ${\mathfrak m}/{\mathfrak m}^2$
respectively. In particular, both $\Tor^S_1(K,K)$ and
$\Tor^R_1(K,K)$ are $r$-dimensional. The proposition follows.
\end{proof}
\begin{proof}[Proof of Theorem~\ref{theo2}]
Let us write $T=K[V]$.
We consider the $T$-module $U$, defined by
$$
U=\{(w_1,w_2,\dots,w_r)\in T[-d_1]\oplus\cdots\oplus T[-d_r]\mid
\sum_{i=1}^r w_if_i=0\}.
$$
Since $I=(f_1,f_2,\dots,f_r)$ is $\tau_G(V)$-regular
(in the sense of Mumford and Castelnuovo), we get
that $U$ is generated in degree $\leq \tau_G(V)+1$.
The module
$$
M=\{(w_1,w_2,\dots,w_r)\in R[-d_1]\oplus\cdots\oplus R[-d_r]\mid
\sum_{i=1}^r w_if_i=0\},
$$
gives an exact sequence
$$
0\to M\to R[-d_1]\oplus\cdots\oplus R[-d_r]\to R\to K\to 0.
$$
We can identify $M/{\mathfrak m}M$ with $\Tor^R_2(K,K)$.
The module $M$ is equal to $U^G$.
We have that $((f_1,f_2,\dots,f_r)U)^G=
(f_1,f_2,\dots,f_r)U^G={\mathfrak m}M$ since $f_1,\dots,f_r$ are
invariant and $G$ is linearly reductive.
We can view
$M/{\mathfrak m}M$ as a submodule of $U/(f_1,\dots,f_r)U$.
It is easy to see that every element of $U$ of degree $\geq 2\tau_G(V)+1$,
must lie in $(f_1,\dots,f_r)U$ since $U$ is generated
in degree $\leq \tau_G(V)+1$ and every polynomial of degree $\geq \tau_G(V)$
lies in $(f_1,\dots,f_r)$.
This shows that 
$$
\deg(\Tor^R_2(K,K))=\deg(M/{\mathfrak m}M)\leq \deg(U/(f_1,\dots,f_r)U)\leq 
2\tau_G(V).
$$
By the previous proposition, we also get
$\deg(\Tor_1^S(R,K))\leq 2\tau_G(V)$.
\end{proof}

\end{document}